\documentclass[11pt]{amsart}
\usepackage{amsmath, amsthm, amssymb}
\usepackage{graphicx}
\usepackage{mathbbol}
\usepackage{txfonts}
\usepackage[usenames]{color}

\newtheorem{thm}{Theorem}[section]
\newcommand{\bt}{\begin{thm}}
\newcommand{\et}{\end{thm}}

\newtheorem{ex}[thm]{Example}

\newtheorem{cor}[thm]{Corollary}   

\newcommand{\bc}{\begin{cor}}
\newcommand{\ec}{\end{cor}}

\newtheorem{lem}[thm]{Lemma}   

\newcommand{\bl}{\begin{lem}}
\newcommand{\el}{\end{lem}}

\newtheorem{prop}[thm]{Proposition}
\newcommand{\bp}{\begin{prop}}
\newcommand{\ep}{\end{prop}}

\newtheorem{defn}[thm]{Definition}
\newtheorem{conj}[thm]{Conjecture}

\newcommand{\bd}{\begin{defn}}    
\newcommand{\ed}{\end{defn}}

\theoremstyle{definition}
\newtheorem{rmrk}[thm]{Remark}   

\newcommand{\br}{\begin{rmrk}}
\newcommand{\er}{\end{rmrk}}

\newtheorem{example}[thm]{Example}

\newcommand{\Ln}{\operatorname{log}}

\newcommand{\be}{\begin{equation}}
\newcommand{\ee}{\end{equation}}

\newcommand{\R}{\mathbb{R}}

\newcommand{\diam}{\operatorname{diam}}

\newcommand{\Lip}{\operatorname{Lip}}

\newcommand{\vol}{\operatorname{Vol}}

 \DeclareMathOperator{\so}{SO}

\def\rr{\mathbb{R}}

\def\isom{\cong}

\def\too{\longrightarrow}

\def\madm{m_{\mathrm{ADM}}}

\def\rsch{r_{\mathrm{Sch}}}

\def\gsch{g_{\mathrm{Sch}}}
\def\gschL{g_{\mathrm{Sch},L}}
\def\msch{m_{\mathrm{Sch}}}
\def\Msch{M_{\mathrm{Sch}}}

\def\RSB{\mathrm{RotSym}^\partial}

\DeclareMathOperator{\dil}{dil}

\begin{document}

\title[Near-equality of Penrose]{Near-equality of the Penrose Inequality for rotationally symmetric Riemannian manifolds}

\author{Dan A. Lee}
\thanks{Lee is partially supported by a PSC CUNY Research Grant and NSF DMS \#0903467.}
\address{CUNY Graduate Center and Queens College}
\email{dan.lee@qc.cuny.edu}

\author{Christina Sormani}
\thanks{Sormani is partially supported by a PSC CUNY Research Grant and NSF DMS \#1006059.}
\address{CUNY Graduate Center and Lehman College}
\email{sormanic@member.ams.org}

\begin{abstract}
This article is the sequel to \cite{Lee-Sormani1}, which dealt with the near-equality case of the Positive Mass Theorem.  We study the near-equality case of the Penrose Inequality for the class of 
complete asymptotically flat rotationally symmetric
Riemannian manifolds with nonnegative scalar curvature whose boundaries are outermost minimal hypersurfaces.  Specifically, we prove that if the Penrose Inequality is sufficiently close to being an equality on one of these manifolds, then it must be close to a Schwarzschild space with an appended cylinder, in the sense of Lipschitz Distance.  Since the Lipschitz Distance bounds the Intrinsic Flat Distance on compact sets, we also obtain a result for Intrinsic Flat Distance, which is a more appropriate distance for more general near-equality results, as discussed in \cite{Lee-Sormani1}.
\end{abstract}

\maketitle

\section{Introduction}\label{introduction}

The (Riemannian) Penrose Inequality states that if $(M^n,g)$ is a complete asymptotically flat manifold of nonnegative scalar curvature whose boundary is an outermost minimal hypersurface, then the ADM mass of $(M,g)$ satisfies
   \be    \madm \ge \frac{1}{2} \left(\frac{|\partial M|}{\omega_{n-1}}\right)^{\frac{n-2}{n-1}},   \ee 
where $|\partial M|$ denotes the (hyper-)area of $\partial M$, and $\omega_{n-1}$ is the (hyper-)area of the standard $(n-1)$-sphere in $\rr^n$.  Furthermore, if equality holds, then $(M,g)$ must be isometric to a Riemannian Schwarzschild manifold (or more precisely, the part lying outside its outermost minimal hypersurface).  The second statement may be thought of as a rigidity theorem, and it is natural to consider the \emph{stability} of this rigidity statement.  That is, if the ratio of the two sides of the inequality is close to one, then
in what sense can we say that the manifold is ``close'' to a Schwarzschild manifold?  

The ADM mass was defined by Arnowitt-Deser-Misner in \cite{ADM-mass}.  The Penrose Inequality is a refinement of the Positive Mass Theorem, which was proven by R.\ Schoen and S.-T.\ Yau in dimensions less than eight  \cite{Schoen-Yau-positive-mass, Schoen-1989}, and by Witten for spin manifolds \cite{Witten-positive-mass, Bartnik-1986}.   The $3$-dimensional Penrose Inequality was first proven by G.\ Huisken and T.\ Ilmanen \cite{Huisken-Ilmanen} (where $|\partial M|$ must be replaced by the area of the largest component of $\partial M$), and later proven in full
by H.\ Bray using a different method \cite{Bray-Penrose}.  Bray and the first author extended Bray's proof to spin manifolds of dimension less than eight \cite{Bray-Lee}.

The problem of stability for the Positive Mass Theorem has been studied by the first author in \cite{Lee-near-equality}, by F.\ Finster with Bray and I.\ Kath in \cite{Bray-Finster, Finster-Kath, Finster},
and by J.\ Corvino in \cite{Corvino}.  Since the general problem of stability is a difficult one, we considered the special case of rotational symmetry (that is, $\so(n)$ symmetry), and in that case the authors were able to prove a comprehensive stability result \cite{Lee-Sormani1}.  
Although there is no
Lipschitz stability in that setting (in the sense of Lipschitz Distance defined below), we formulated the stability
in terms of Intrinsic Flat Distance, a notion defined by the
second author and S.\ Wenger in \cite{SorWen2}.  This article is a natural extension of that investigation to the problem of stability for the Penrose Inequality.

We find that stability fails for the Penrose Inequality
in the sense that manifolds with almost equality in the
Penrose Inequality are not necessarily close to Schwarzschild
space in any reasonable topology.  Specifically, they may instead be close to Schwarzschild spaces with a cylinder of arbitrary length appended to the boundary. 
See Example~\ref{ex-penrose-well}  depicted in 
Figure~\ref{fig-thm-penrose} below.   
This behavior is to be expected because the rigidity in the Penrose Inequality only applies to manifolds whose boundaries are outermost minimal hypersurfaces.  Since this class of Riemannian manifolds is not closed in any reasonable topology, one does not expect a true stability result.  
Before we state our result, we provide a few definitions.

\begin{defn}  \label{def-rot-sym}
Given $n\ge 3$,
let $\RSB_n$ be the class of complete $n$-dimensional 
rotationally symmetric (that is, $\so(n)$ symmetry) smooth complete Riemannian manifolds
of nonnegative scalar curvature whose connected nonempty boundaries are outermost (and outer-minimizing) minimal hypersurfaces. \end{defn}

The class $\RSB_n$ includes the Schwarzschild spaces.  The nonnegative scalar curvature condition corresponds to the physical assumption of nonnegative mass density in the time-symmetric setting.
The outermost condition is included here, just
as it is in the Penrose Inequality, because
complicated geometry can ``hide'' behind a minimal hypersurface
without affecting the ADM mass, \emph{cf.\ }\cite{Gibbons, Huisken-Ilmanen}).  The outermost minimal boundary is often called an apparent horizon.  
Note that we need not explicitly assume asymptotic flatness here 
because finite ADM mass in $\RSB_n$ implies asymptotic flatness.

Recall the Lipschitz Distance between metric spaces as defined in \cite{Gromov-metric}:
\begin{defn}
Given two metric spaces $(X,d_X)$ and $(Y,d_Y)$, the Lipschitz Distance
\be\label{eqn-defn-lip}
d_{\mathcal{L}}(X,Y):=\inf\left\{ |\Ln \dil (\varphi)| + |\Ln \dil (\varphi^{-1})|\,:\, \varphi:X\too Y\text{\emph{ is bi-Lipschitz}}\right\}
\ee
where 
\be
\dil (\varphi)=\sup_{x_i\in X} \frac{ d_Y(\varphi(x_1),\varphi(x_2))}{d_X(x_1,x_2)}.
\ee
\end{defn}

\begin{defn}
For any $m>0$, define $\Msch^n(m)$ to be the part of the $n$-dimensional Schwarzschild space of mass $m$ that lies outside its outermost minimal hypersurface.  For any $L\in[0,\infty)$, define $\Msch^n(m,L)$ to be $\Msch^n(m)$ with the cylinder $[-L,0]\times S^{n-1}((2m)^{\frac{1}{n-2}})$ glued to its outermost minimal hypersurface boundary, where $S^{n-1}((2m)^{\frac{1}{n-2}})\isom\partial \Msch^n(m)$ denotes the standard $(n-1)$-sphere of radius $(2m)^{\frac{1}{n-2}}$. 
\end{defn}

\begin{thm} \label{thm-penrose}
Let $n\ge3$.  For any $\epsilon>0$, there exists $\delta>0$ such that if $M\in \RSB_n$ satisfies
   \be   
    \madm \le \frac{1+\delta}{2} \left(\frac{|\partial M|}{\omega_{n-1}}\right)^{\frac{n-2}{n-1}},   
    \ee 
then
   \be   
    d_{\mathcal{L}}\left(M, \Msch^n (m_0, L)\right) < \epsilon,   
    \ee 
where $m_0= \frac{1}{2} \left(\frac{|\partial M|}{\omega_{n-1}}\right)^{\frac{n-2}{n-1}}$, and $L$ is the \emph{depth} of $M$, as defined in Section \ref{depth}.
\end{thm} 

In other words, if the ratio of the two sides of the Penrose Inequality is close enough to one, then the space must be close to an appended Schwarzschild space in Lipschitz Distance.  See Theorem~\ref{thm-precise} for a more precise statement implicitly describing the dependence of $\delta$
on $\epsilon$.

As alluded to above, the reason why we need the appended Schwarzschild spaces is that we may have a sequence of manifolds in $\RSB_n$ whose limit has a boundary that is not outermost.  This situation is described in Example~\ref{ex-penrose-well}  and depicted in Figure~\ref{fig-thm-penrose} below.   Example~\ref{ex-not-C2} demonstrates that
one cannot improve the Lipschitz convergence to $C^2$ convergence.
 
\begin{figure}[h] 
   \centering
   \includegraphics[width=5in]{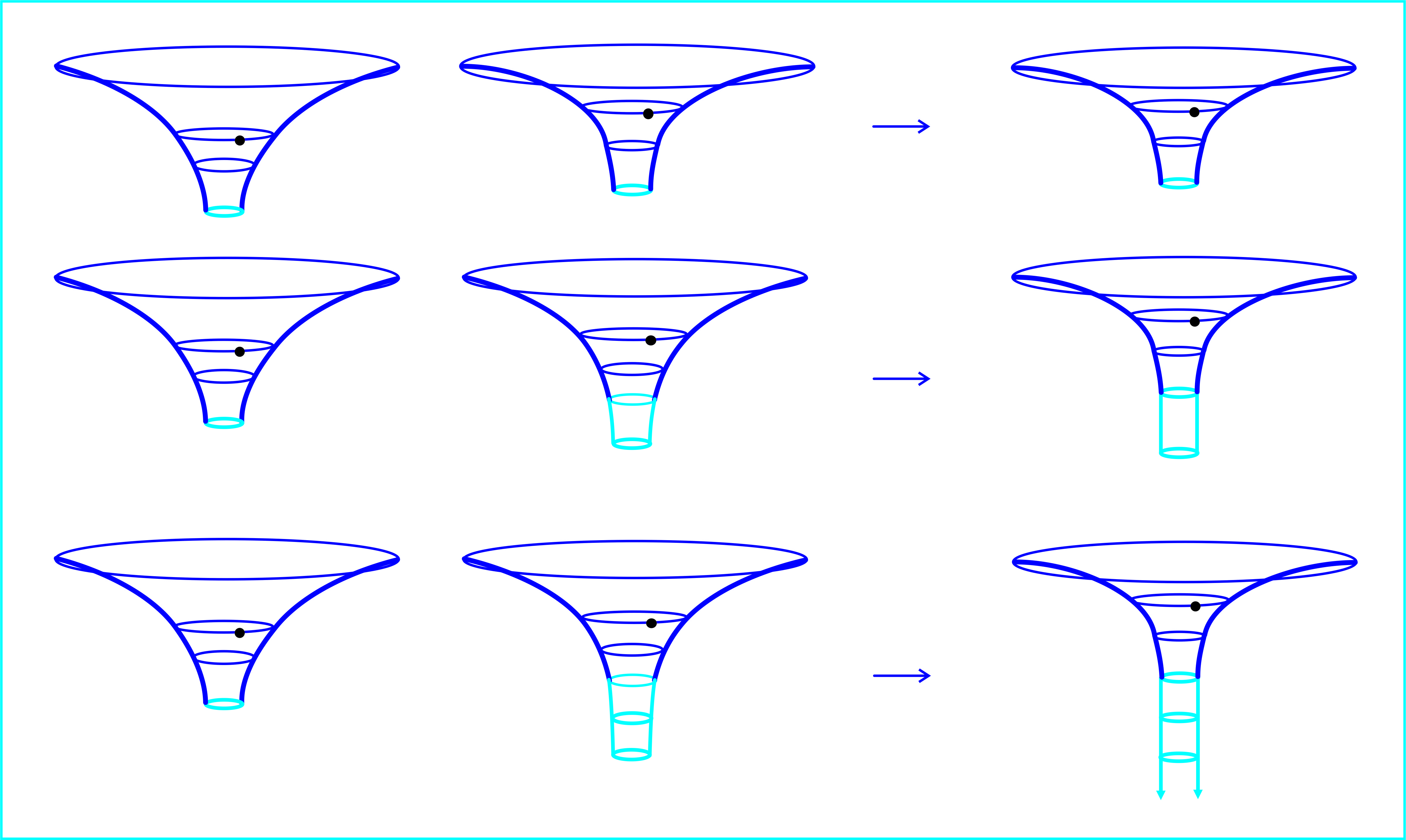} 
   \caption{Sequences approaching various limits, $\Msch(m,L)$.
   }
   \label{fig-thm-penrose}
\end{figure}

Note that Theorem \ref{thm-penrose} does not hold without the assumption of rotational symmetry, as seen in Example~\ref{ex-many-wells}.  The basic reason why we obtain such a strong result in Theorem \ref{thm-penrose} is that rotational symmetry combined with the existence of a boundary of fixed size rules out the possibility of a thin, deep gravity well.  These deep gravity wells are the reason why neither Lipschitz Distance nor even Gromov-Hausdorff Distance provide a useful topology for a more general theorem.  This is discussed in our earlier work \cite{Lee-Sormani1}, in which we propose Intrinsic Flat Distance as an appropriate topology.  

The Intrinsic Flat Distance was defined by the second 
author and S.\ Wenger in \cite{SorWen2} based upon
work of L.\ Ambrosio and B.\ Kirchheim \cite{AK}.  It estimates the distances
between compact Riemannian $n$-manifolds by filling in (in the sense of metric space isometric embeddings) the space between
them with a countably $\mathcal{H}^{n+1}$ rectifiable metric space and measuring 
the $\mathcal{H}^{n+1}$ measure of the filling space and the $\mathcal{H}^{n}$ measure of any excess boundary.
It was proven in \cite{SorWen2} that the Intrinsic Flat Distance between
two Riemannian manifolds with boundary can be bounded in terms of the
Lipschitz Distance between them, their diameters, their
volumes and the (hyper-)areas of their boundaries.  In
recent of work of the second author and S.\ Lakzian, an
explicit filling manifold is constructed between two
spaces which are close in the Lipschitz sense \cite{Lakzian-Sormani}.
See Section \ref{intrinsic-flat} for the precise statement.
Using this bound we conclude the following:

\begin{thm}\label{thm-penrose-IF}
Let $n\ge3$, and let $0<A_0<A_1$.   For any $D>0$ and $\epsilon>0$, there exists $\delta>0$ such that if 
$M\in\RSB_n$ satisfies $|\partial M|=A_0$ and
   \be \label{aaaa}
      \madm\le \frac{1+\delta}{2} \left(\frac{A_0}{\omega_{n-1}}\right)^{\frac{n-2}{n-1}}, 
        \ee 
then
   \be   \label{bbbb}
    d_{\mathcal{F}}\left(T_D(\Sigma_1)\subset M, T_D(\bar{\Sigma}_1)\subset \Msch^n (m_0, L)\right) < \epsilon,  
     \ee 
where $m_0= \frac{1}{2} \left(\frac{|\partial M|}{\omega_{n-1}}\right)^{\frac{n-2}{n-1}}$, $L$ is the depth of $M$, 
$\Sigma_1$ and $\bar{\Sigma}$ are the respective symmetric spheres of area $A_1$ in $M$ and $\Msch(m_0,L)$, $T_D$ denotes the tubular neighborhood of radius $D$, and  $d_{\mathcal{F}}$ denotes Intrinsic Flat Distance.
\end{thm}

For more discussion of Intrinsic Flat Distance and its application to the study of nonnegative scalar curvature, see \cite{Lee-Sormani1}.  Note that a scale-invariant version of this 
result would require a scalable version of Intrinsic Flat 
Distance.  This is being developed by J.\ Basilio \cite{Basilio1}.

In the final section of the paper we provide examples 
and
propose that some version of Theorem~\ref{thm-penrose-IF} holds without
the restriction of rotational symmetry.  See Conjecture~\ref{conj-penrose-IF}
and subsequent remarks.

The authors would like to thank Lars Andersson, Hubert Bray, 
Piotr Chru\'{s}ciel, Lan-Hsuan Huang, Gerhard Huisken, Tom Ilmanen, 
James Isenberg, John Lee, and Shing-Tung Yau for their interest in this work.   


\section{Basic facts about $\RSB_n$}\label{background}

\subsection{Geodesic coordinates}
In this article we consider Riemannian manifolds $(M,g)$ in $\RSB_n$,
defined in Definition~\ref{def-rot-sym}.
Since $(M,g)$ is rotationally symmetric 
we can write its metric in
geodesic coordinates, as $g=ds^2 + r(s)^2 g_0$ for some function 
$r:[0,\infty) \to (0,\infty)$,
where $g_0$ is the standard metric on the $(n-1)$-sphere
and $s$ is the distance from the boundary, 
$\partial M$.   

Let $\Sigma$ be the symmetric sphere that is a distance
$s$ from the boundary.
We then have the following formulae for the (hyper-)area and mean curvature:
\begin{alignat}{1}
|\Sigma|&=\omega_{n-1} r^{n-1}\\
H_\Sigma & = \frac{n-1}{r}\frac{dr}{ds}.
\end{alignat}
At a point $p\in \Sigma$ the scalar curvature is
\be \label{eqn-18}
\mathrm{R} = \frac{n-1}{r^2}\left((n-2)\left[1- \left(\frac{dr}{ds}\right)^2\right] - 2r\frac{d^2r}{ds^2}\right).
\ee

Since $\partial M$ is an outermost minimal hypersurface,  $\left.\frac{dr}{ds}\right|_{s=0}=0$, and $\frac{dr}{ds}$ is either positive for all $s>0$ or negative for all $s>0$.  Since $\partial M$ is outer-minimizing, it must be the case that
\be\label{eqn-f'>0}
\frac{dr}{ds}>0 \qquad \forall s\in (0,\infty).  
\ee

Recall the definition of the Hawking mass of a surface $\Sigma$
in a three-dimensional manifold:
\be
m_H(\Sigma) := \frac{1}{2}\sqrt{\frac{|\Sigma|}{\omega_{2}}}\left(1-\frac{1}{\omega_{2}}\int_\Sigma \left(\frac{H}{2}\right)^{2}\right).
 \ee
We define a natural Hawking mass function on symmetric spheres $\Sigma$ in $M\in \RSB_n$ that agrees with the usual definition in dimension three.
\be \label{define-hawking}
m_H(\Sigma) := \frac{r^{n-2}}{2}\left(1-\left(\frac{dr}{ds}\right)^2\right). 
\ee
Alternatively, we may view $m_H$ as a function of $s$.  Applying (\ref{eqn-18}), one can compute
\be \label{eqn-hawking-2}
\frac{dm_H}{ds} = \frac{r^{n-1}}{2(n-1)}\frac{dr}{ds} \mathrm{R}.
\ee
Since we are studying manifolds with $\frac{dr}{ds}>0$ for $s>0$, we observe that the monotonicity of the Hawking mass, usually called \emph{Geroch monotonicity}:
\be\label{geroch}
\frac{dm_H}{ds}\ge 0,
\ee
is actually \emph{equivalent} to nonnegativity of $\mathrm{R}$.

The ADM mass of $M$ is defined as the limit of the Hawking mass function:
\be
\madm(M) :=\lim_{s\to \infty} m_H  \in [0,\infty].
\ee
For $M\in\RSB_n$, this agrees with the usual definition of the ADM mass.  Note that the Penrose Inequality for manifolds in $\RSB_n$ then follows immediately from Geroch monotonicity \eqref{geroch}.

\subsection{Graphical coordinates} \label{graphical-coordinates}
Since  $\frac{dr}{ds}>0$ for $s>0$, the map $s\mapsto r(s)$ defines a smooth change of coordinates away from $\partial M$.  In $r$ coordinates, the metric takes the form
   \be \label{g-radial}  g = \left(\frac{ds}{dr}\right)^2 dr^2 + r^2 g_0.   \ee 
By Geroch monotonicity \eqref{geroch}, we see that $\frac{dr}{ds}<1$ everywhere.  We now choose $z(r)$ to be an increasing function (determined up to a constant) such that
   \be \label{define-z}    1+  \left(\frac{dz}{dr}\right)^2 = \left(\frac{ds}{dr}\right)^2,   \ee 
so that
   \be   g =  \left(1+  \left(\frac{dz}{dr}\right)^2\right) dr^2 + r^2 g_0.   \ee 
Note that this formula exhibits a Riemannian isometric embedding of $(M,g)$ into $\rr^{n+1}$ as the graph of the radial function $z(r)$.  That is, we may view
\be
M = \left\{ (x',x^{n+1})\in\rr^{n+1}\,:\, x^{n+1} = z (|x'|) \right\}.
\ee
For this reason, we call this choice of coordinates \emph{radial graphical coordinates}.

We also observe that $r\mapsto z(r)$ defines a smooth change of coordinates away from $\partial M$, and in $z$ coordinates, the metric is
\be
\label{g-vertical}
g =  \left(\frac{ds}{dz}\right)^2 dz^2 + r(z)^2 g_0. 
\ee
Also note that unlike the $r$-coordinates, these $z$-coordinates are nonsingular at $\partial M$.  We call this choice of coordinates \emph{vertical graphical coordinates}.

\subsection{Appended Schwarzchild spaces}\label{schwarzschild}
We recall that the $n$-dimensional Schwarzschild space $\Msch^n(m)$ of mass $m$ may be described in radial graphical coordinates as the manifold $[(2m)^{\frac{1}{n-2}},\infty)\times S^{n-1}$ with the metric
   \be\label{SC-radial}    \gsch =  \left(1 -\frac{2m}{r^{n-2} }\right)^{-1}dr^2   + r^2 g_0.    \ee 
In vertical graphical coordinates, the metric becomes
   \be\label{SC-vertical}    \gsch =  \frac{r(z)^{n-2}}{2m} dz^2   + r(z)^2 g_0,    \ee 
where 
\be\label{SC-define-z}
\frac{dz}{dr}=\sqrt{\frac{2m}{r^{n-2}-2m}}.
\ee
  By convention, we choose vertical graphical coordinates so that $\partial\Msch(m)$ has \hbox{$z=0$}.  For $L\ge0$, the appended Schwarzschild space $\Msch(m,L)$ is obtained by gluing $\Msch(m)$ to the space $[-L,0]\times S^{n-1}$ with the cylindrical metric
   \be \label{cylinder}   g_{\mathrm{cyl}} = dz^2 + (2m)^{\frac{2}{n-2}}g_0   \ee 
along their common sphere at $z=0$.


\section{Lipschitz Estimates}\label{lipschitz}

In this section we define the depth
and prove all Lipschitz estimates needed to prove our main results.

\subsection{Depth} \label{depth}

Let $M\in \RSB_n$.  For convenience,
we express   
the important quantities $\madm$ and $|\partial M|$ in different units, defining
   \begin{align}
   m_0&:=\frac{1}{2}\left(\frac{|\partial M|}{\omega_{n-1}}\right)^{\frac{n-2}{n-1}}>0\\
      r_0&:= r(0)=(2m_0)^{\frac{1}{n-2}} = \left(\frac{|\partial M|}{\omega_{n-1}}\right)^{\frac{1}{n-1}}>0 \label{r_0}\\
       r_1&:= (2\madm)^{\frac{1}{n-2}}>0.\label{r_1}  
       \end{align}
These quantities will be applied to define the depth and again in subsequent
proofs.
       
 Define $\delta$ so that
   \be  \label{assumption}
   \madm = \frac{1+\delta}{2} \left(\frac{|\partial M|}{\omega_{n-1}}\right)^{\frac{n-2}{n-1}},  
   \ee 
Using the definitions above, we can rewrite this as
   \be    r_1^{n-2} = (1+\delta) r_0^{n-2}.   \ee 
 We now define $r_\delta$ and $A_\delta$ so that
   \begin{align}   r_\delta&:= (1+\sqrt{\delta})^{\frac{1}{n-2}}r_0\\
   A_\delta&:=\omega_{n-1}r_\delta^{n-1}.
    \end{align}

Recall that $\Msch(m_0)$ may be isometrically embedded into
$\rr^{n+1}$ with its boundary embedded at the level $z=0$
and radius $r=r_0$, as described in Section~\ref{schwarzschild}.
We can also isometrically embed $M$ into $\rr^{n+1}$ via the function $z(r)$ introduced in Section \ref{graphical-coordinates}, and this embedding is defined only up to constant.  Let us choose this constant so that $M$ intersects the standard embedding of $\Msch(m_0)$ in $\rr^{n+1}$ at $r=r_\delta$.   See Figure~\ref{LS2-fig2} below.

\begin{figure}[h] 
   \centering
   \includegraphics[width=4in]{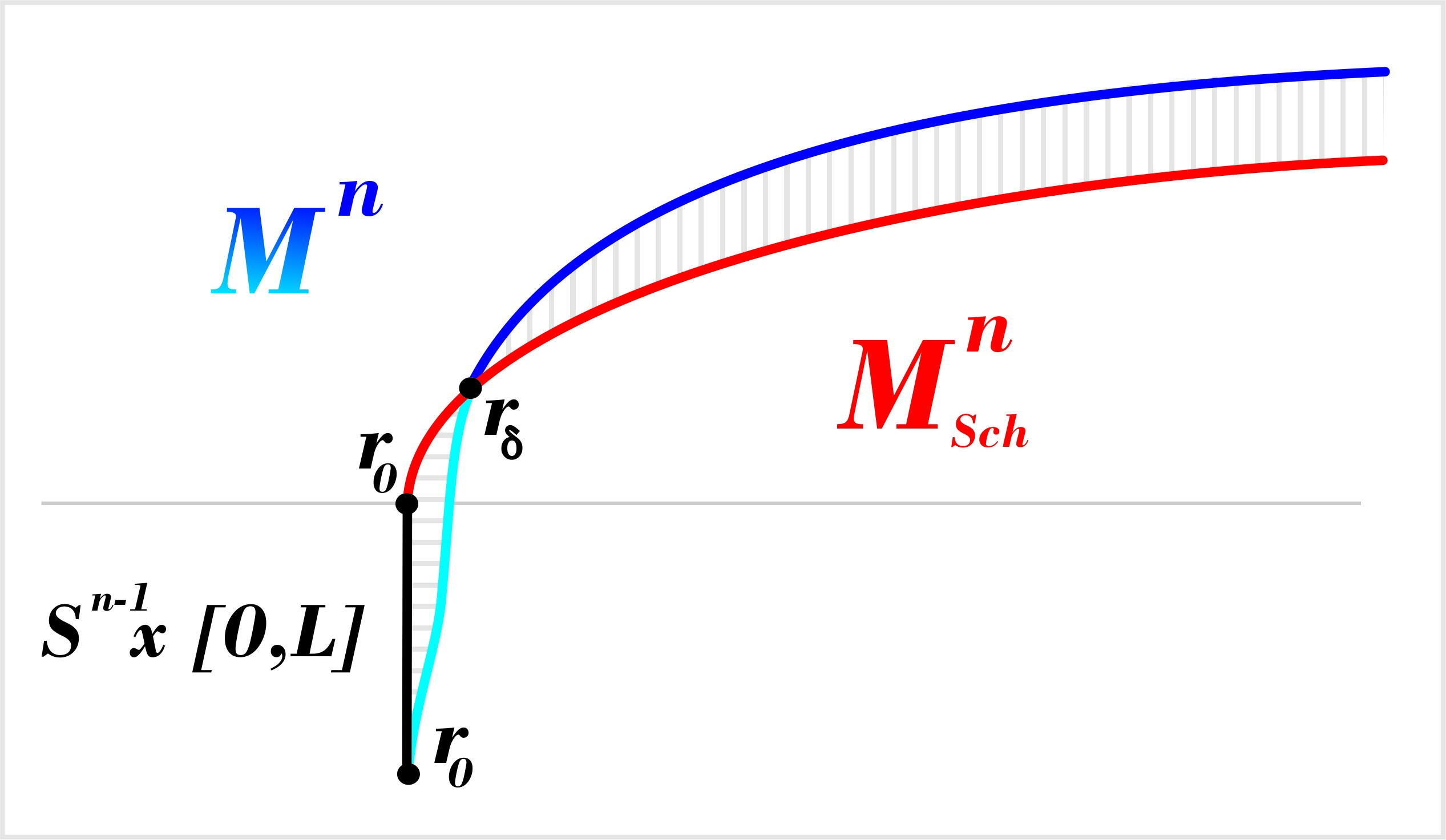} 
   \caption{The depth $L$ and the bi-Lipschitz map $\varphi$ are defined by embedding $M$ and $\Msch(m_0)$ into $\rr^{n+1}$ so that they intersect at $r=r_\delta$.
   }
   \label{LS2-fig2}
\end{figure}

\begin{defn}
We define the {\em depth} of $M\in \RSB_n$ 
to be $L$ such that when $M$ is embedded into $\rr^{n+1}$
as described above and depicted in Figure~\ref{LS2-fig2},
$\partial M$ embeds into $r=r_0$ and $z=-L$. 
\end{defn}

\begin{rmrk}\label{rmrk-depth}
We note that although the depth is an invariant of $M\in\RSB_n$ 
which could plausibly be generalized to  non-rotationally symmetric 
spaces, the definition is not quite natural in the sense that it depends 
on an arbitrary choice for what $r_\delta$ should be.  As one can see from the proof in the next section, if we were to replace 
$\sqrt{\delta}$ by $\delta^{1/3}$ in the definition of depth, 
Theorem~\ref{thm-penrose} would still hold true.  Nevertheless, 
the concept of depth may be worthy of further investigation.
\end{rmrk}

\begin{lem}
The depth $L$ is always nonnegative.  
\end{lem}

\begin{proof}
Consider the functions $r(s)$ and $z(r)$ as defined for the manifold $M$.  
By Geroch monotonicity (\ref{geroch}), we have 
   \be    m_H(\partial M) \le m_H(\Sigma)  \le \madm   \ee 
for any symmetric sphere $\Sigma$.  By the definitions of Hawking mass \eqref{define-hawking} and $r_0$ and $r_1$ in \eqref{r_0} and \eqref{r_1}, this becomes
   \be\label{pinching1}    r_0^{n-2} \le r^{n-2}\left(1-\left(\frac{dr}{ds}\right)^2\right)\le r_1^{n-2}.   \ee 
Therefore
\begin{equation}\label{pinching2}
1- \left(\frac{r_0}{r}\right)^{n-2}\ge
\left(\frac{dr}{ds}\right)^2 \ge 1- \left(\frac{r_1}{r}\right)^{n-2}.
\end{equation}
By combining \eqref{pinching2} and \eqref{define-z}, we see that for all $r\ge r_0$, 
\begin{gather}
1- \left(\frac{r_0}{r}\right)^{n-2} \ge \left(1+\left(\frac{dz}{dr}\right)^2\right)^{-1}\\
\frac{r^{n-2}}{r^{n-2}- r_0^{n-2}}\le  1+\left(\frac{dz}{dr}\right)^2 \\
\frac{dz}{dr} \ge \sqrt{ \frac{r_0^{n-2}}{r^{n-2}- r_0^{n-2}}}= \sqrt{ \frac{2m_0}{r^{n-2}- 2m_0}}.
\end{gather}
Using \eqref{SC-define-z} and the fact that the graphs of $M$ and $\Msch(m_0)$ intersect at $r=r_\delta$, we can integrate the above inequality to see that $M$ cannot lie above $\Msch(m_0)$ in $\rr^{n+1}$ in the region where $r\le r_\delta$.  In particular, this means that $\partial M$ cannot lie above $\partial\Msch(m_0)$, so that $-L=z(r_0)\le 0$, as desired.
\end{proof}

\subsection{Proof of Theorem~\ref{thm-penrose}}

Theorem~\ref{thm-penrose} is an easy consequence
of the following theorem.

\begin{thm}\label{thm-precise}
Let $n\ge3$, and let $M\in \RSB_n$.  Define $\delta$ so that
   \be   
   \madm = \frac{1+\delta}{2} \left(\frac{|\partial M|}{\omega_{n-1}}\right)^{\frac{n-2}{n-1}},  
   \ee 
and assume $0\le \delta<1$.  Define $m_0= \frac{1}{2} \left(\frac{|\partial M|}{\omega_{n-1}}\right)^{\frac{n-2}{n-1}}$, $A_\delta=(1+\sqrt{\delta})^{\frac{n-1}{n-2}}|\partial M|$, and let $L$ be depth of $M$.  Then there exists a rotationally symmetric
bi-Lipschitz map $\varphi: M \to \Msch(m_0,L)$ 
that maps symmetric spheres of area
$A\ge A_\delta$ to symmetric spheres of area $A$
such that for any tangent vector $v\in TM$,
\be \label{g-ratio}
h_\delta^{-1} \le \frac{g(v,v)}{\gschL(\varphi_* v, \varphi_*v)}
\le h_\delta
\ee
where
\be
h_\delta=\max\left\{(1+\sqrt{\delta})^{\frac{2}{n-2}},\, (1+\sqrt{\delta})(1+\delta),\,
(1-\sqrt{\delta})^{-1} \right\}.
\ee
\end{thm}

Theorem~\ref{thm-penrose} follows from this theorem
because  $\sqrt{h_\delta}$ is then a bound on $\Lip(\varphi)$
and $\Lip(\varphi^{-1})$, so that 
\be
d_{\mathcal{L}}(M, \Msch(m_0,L) ) \le 2 \log \sqrt{h_\delta}=\log(h_\delta),
\ee
which converges to $0$ as $\delta\to0$.

\begin{proof}
We will prove Theorem \ref{thm-precise} by explicitly constructing the bi-Lipschitz map.  Let $(M,g)$, $\delta$, $m_0$, $A_\delta$, and $L$ be as described in the statement of the theorem.  We define $r_0$, $r_1$, and $r_\delta$ as in Section \ref{depth}, and we isometrically embed $M$ and $\Msch(m_0)$ into $\rr^{n+1}$ as in Section \ref{depth}. 
Then we embed $\Msch(m_0,L)$ so that $\partial M=\partial\Msch(m_0,L)$ and such that the Schwarzschild domain within $\Msch(m_0,L)$ 
agrees with the embedding of $\Msch(m_0)$.  See Figure~\ref{LS2-fig2} above.

We define the bi-Lipschitz map $\varphi: M\too \Msch(m_0,L)$ as follows.  In the region where $r\ge r_\delta$, $\varphi$ projects $M$ vertically in $\rr^{n+1}$ to $\Msch(m_0, L)$.  That is, under this map, the two manifolds have the same $r$ coordinates and sphere coordinates, but different $z$ coordinates.
In particular $\varphi$ preserves the areas of symmetric
spheres of area $A\ge A_\delta$.
In the region where $r \le r_\delta$, $\varphi$ projects $M$ horizontally in $\rr^{n+1}$ to $\Msch(m_0, L)$.  That is, under this map, the two manifolds have the $z$ coordinates and sphere coordinates, but different $r$ coordinates.  The vertical and horizontal lines in Figure~\ref{LS2-fig2} depict the map $\varphi$.  Since the two definitions match up at $r=r_\delta$, this map is clearly bi-Lipschitz.  

We now prove (\ref{g-ratio}) for the region where $r\ge r_\delta$.  
Recall that in this region, $\varphi$ preserves both the $r$ coordinate and the sphere coordiantes, so it suffices to assume $v$ is purely radial.  Note that the assumption that $\delta<1$ implies that $r_0\le r_1<r_\delta$, so that when $r\ge r_\delta$, \eqref{pinching2} tells us that 
\be\label{pinching2again}
1- \left(\frac{r_0}{r}\right)^{n-2}\ge
\left(\frac{dr}{ds}\right)^2 \ge 1- \left(\frac{r_1}{r}\right)^{n-2}>0.
\ee
By comparing \eqref{g-radial} to \eqref{SC-radial}, the inequalities in \eqref{pinching2again} tell us that for any radial 
tangent vector $v\in TM$ based at a point with $r\ge r_\delta$,  
\be
 \frac{1- \frac{2m_0}{r^{n-2}}}{1- \left(\frac{r_0}{r}\right)^{n-2}}  \le \frac{g(v,v)}{\gschL(\varphi_* v, \varphi_*v)} \le \frac{1- \frac{2m_0}{r^{n-2}}} {1- \left(\frac{r_1}{r}\right)^{n-2}}.  
 \ee
 By the definition of $r_0$, we have
  \be
   1  \le \frac{g(v,v)}{\gschL(\varphi_* v, \varphi_*v)} \le \frac{1- \left(\frac{r_0}{r}\right)^{n-2}} {1- \left(\frac{r_1}{r}\right)^{n-2}}.
   \ee
One can see that over $r\in[r_\delta,\infty)$, the right hand side is maximized at $r=r_\delta$, so that 
   \be    \frac{g(v,v)}{\gschL(\varphi_* v, \varphi_*v)} \le \frac{r_\delta^{n-2} - r_0^{n-2}}{r_\delta^{n-2}-r_1^{n-2}}= \frac{1}{1-\sqrt{\delta}}.   
   \ee

We now prove (\ref{g-ratio}) for the region
where $r_0\le r\le r_\delta$.  Since in this region, $\varphi$ changes the $r$-coordinate, we must separately consider vectors that are tangent to the symmetric spheres, as well as ones that are orthogonal to the symmetric spheres.  However, the image of a symmetric sphere in $M$ with radius $r_0\le r\le r_\delta$ under $\varphi$ is another symmetric sphere with radius $r_0\le r\le r_\delta$.  Therefore it is clear that for any $v\in TM$ that is tangent to a symmetric sphere and based at a point with $r_0\le  r\le r_\delta$,
\begin{gather}
 \left(\frac{r_0}{r_\delta}\right)^2\le \frac{g(v,v)}{\gschL(\varphi_* v, \varphi_*v)}\le \left(\frac{r_\delta}{r_0}\right)^2\\
(1+\sqrt{\delta})^{-\frac{2}{n-2}} \le\frac{g(v,v)}{\gschL(\varphi_* v, \varphi_*v)}\le (1+\sqrt{\delta})^{\frac{2}{n-2}}.
\end{gather}

Now let us consider tangent vectors that are radial, that is, orthogonal to the symmetric spheres.
Using \eqref{pinching1}, we see that for any $r\ge r_0$,
   \be \left(\frac{r_0}{r}\right)^{n-2} \le \left(1-\left(\frac{dr}{ds}\right)^2\right)\le \left(\frac{r_1}{r}\right)^{n-2}.   \ee 
By \eqref{define-z}, this may be written 
  \be \label{pinching3}
  \left(\frac{r_0}{r}\right)^{n-2} \le \left(\frac{dz}{ds}\right)^2\le \left(\frac{r_1}{r}\right)^{n-2}.   \ee 
Consider a radial tangent vector $v\in TM$ based at a point with $r_0\le r\le r_\delta$.  We consider two cases.  
In the first case, $z(r)\ge0$, so that $\varphi$ projects the point to a point on $\Msch(m_0)$.  Recall that in this region, $\varphi$ does not change $z$.  So by comparing \eqref{g-vertical} to \eqref{SC-vertical},  the inequalities in \eqref{pinching3} tell us that
\begin{gather}
   \left(\frac{r(z)}{r_0}\right)^{n-2} \frac{2m_0}{\rsch(z)^{n-2}}  \ge  \frac{g(v,v)}{\gschL(\varphi_* v, \varphi_*v)} \ge 
 \left(\frac{r(z)}{r_1}\right)^{n-2} \frac{2m_0}{\rsch(z)^{n-2}} \\
   \left(\frac{r_\delta}{r_0}\right)^{n-2}  \ge  \frac{g(v,v)}{\gschL(\varphi_* v, \varphi_*v)} \ge 
 \left(\frac{r_0}{r_1}\right)^{n-2}  \left(\frac{r_0}{r_\delta}\right)^{n-2}\\
1+\sqrt{\delta} \ge  \frac{g(v,v)}{\gschL(\varphi_* v, \varphi_*v)} \ge (1+\delta)^{-1}(1+\sqrt{\delta})^{-1},
\end{gather}
 where $\rsch(z)$ describes how the $r$ coordinate on $\Msch(m_0)$ depends on $z$.  In the second case,  $v$ is based at a point with $z(r)<0$, so that it projects to the cylindrical part of $\Msch(m_0, L)$.   By comparing \eqref{g-vertical} to \eqref{cylinder}, then inequalities in \eqref{pinching3} tell us that
   \begin{gather}  \label{cyl-part}  \left(\frac{r(z)}{r_0}\right)^{n-2} \ge \frac{g(v,v)}{\gschL(\varphi_* v, \varphi_*v)}\ge \left(\frac{r(z)}{r_1}\right)^{n-2}\\
 1+\sqrt{\delta} \ge \frac{g(v,v)}{\gschL(\varphi_* v, \varphi_*v)}\ge (1+\delta)^{-1}.
   \end{gather} 
So we have 
(\ref{g-ratio}) for vectors tangent to symmetric spheres
or perpendicular to them. In this situation, one can also see that
it holds for their linear combinations.
\end{proof}

\subsection{Lipschitz estimates for tubular neighborhoods}

The following proposition is useful for the proof of Theorem \ref{thm-penrose-IF}.

\begin{prop}\label{Lip-Tube}
Let $n\ge3$, and let $0<A_0<A_1$.   For any $D>0$ and $\epsilon>0$, there exists $\delta>0$ such that for all 
$M\in\RSB_n$ satisfying $|\partial M|=A_0$ and
   \be 
      \madm< \frac{1+\delta}{2} \left(\frac{A_0}{\omega_{n-1}}\right)^{\frac{n-2}{n-1}}, 
        \ee 
we have
   \be   \label{claim}
    d_{\mathcal{L}}\left(T_D(\Sigma_1)\subset M, T_D(\bar{\Sigma}_1)\subset \Msch^n (m_0, L)\right) < \epsilon,  
     \ee 
where $m_0= \frac{1}{2} \left(\frac{|\partial M|}{\omega_{n-1}}\right)^{\frac{n-2}{n-1}}$, $L$ is the depth of $M$, $\Sigma_1$ and $\bar{\Sigma}$ are the respective symmetric spheres of area $A_1$ in $M$ and $\Msch(m_0,L)$, $T_D$ denotes the tubular neighborhood of radius $D$, and  $d_{\mathcal{L}}$ denotes Lipschitz Distance.
\end{prop}

\begin{proof}
Following the proof of Theorem \ref{thm-penrose}, for small enough $\delta$ there exist
 $L\ge0$ and a bi-Lipschitz map $\varphi:M\too \Msch(m_0,L)$ such that 
   \be    d_{\mathcal{L}}\left(T_D(\Sigma_1) , \varphi(T_D(\Sigma_1)) \right) <\epsilon/2.    \ee 
While $\varphi(\Sigma_1)=\bar{\Sigma}_1$, it is not true that $\varphi(T_D(\Sigma_1))=T_D(\bar{\Sigma}_1)$.  
We can define a bi-Lipschitz map $\psi:\varphi(T_D(\Sigma_1))\too T_D(\bar{\Sigma}_1)$ that preserves $\bar{\Sigma}_1$ while scaling the distance to $\bar{\Sigma}_1$ by a constant factor on the part inside $\bar{\Sigma}_1$ and another constant factor on the outside part.  Since we know that the ratios between distances in $M$ and $\Msch(m_0,L)$ tend to $1$ as $\delta\to0$, we can arrange for these constant factors to be close to $1$.  This controls $\gschL(v,v)/\gschL(\psi_* v,\psi_* v)$ for radial tangent vectors $v$.  And since $\frac{d\rsch}{ds}\le1$, we can also control $\gschL(v,v)/\gschL(\psi_* v,\psi_* v)$ for $v$ that are tangent to the symmetric spheres.  
So we can choose $\delta$ small enough so that
   \be   d_{\mathcal{L}}(\varphi(T_D(\Sigma_1)),T_D(\bar{\Sigma}_1))<\epsilon/2,   \ee 
completing the proof of \eqref{claim}.
\end{proof}


\section{Intrinsic Flat Distance}

\subsection{Review of the Intrinsic Flat Distance}\label{intrinsic-flat}

The Intrinsic Flat Distance measures the distances between 
compact oriented Riemannian manifolds by filling in the space 
between them.  This notion was first defined in work of the second author
and S.\ Wenger in \cite{SorWen2} applying the theory of integral currents on metric
spaces developed by Ambrosio-Kirchheim in \cite{AK}.   However,
here we estimate the intrinsic flat distance using only an understanding
of Riemannian geometry.

Given two compact oriented Riemannian manifolds $M^n_1$ and $M^n_2$ with boundary, and metric isometric embeddings $\psi_i: M_i \to Z$ into some Riemannian
manifold (possibly piecewise smooth with corners)
an upper bound for the Intrinsic Flat Distance is attained as follows:
\be \label{eqn-def-intrinsic-flat-1}
d_{\mathcal{F}}(M^n_1, M^n_2) \le \vol_{n+1}(B^{n+1}) +\vol_n(A^n)
\ee
where $B^{n+1}$ is an oriented region in $Z$ and $A^n$ is
defined so that the oriented integrals satisfy
\be \label{eqn-def-intrinsic-flat-2}
\int_{\psi_1(M_1)}\omega -\int_{\psi_2(M_2)}\omega=\int_{\partial B}\omega + \int_A\omega 
\ee
for any differential $n$-form $\omega$ on $Z$.   We call
$B^{n+1}$ a {\em filling manifold} between $M_1$ and $M_2$
and $A^n$ the {\em excess boundary}.

Recall that a metric isometric embedding, 
$\psi: M \to Z$ is a map such that
\be\label{eqn-def-metric-isom-embed}
d_Z(\psi(x), \psi(y)) = d_M(x,y) \qquad \forall x,y \in M.
\ee
This is significantly stronger than a Riemannian isometric embedding
which preserves only the Riemannian structure and thus lengths
of curves but not distances between points as in 
(\ref{eqn-def-metric-isom-embed}).   For example, $S^1$ has
a metric isometric embedding
as the equator of a hemisphere but only a Riemannian isometric 
embedding as a circle in a plane. 

Theorem 5.6 of \cite{SorWen2} states that
the Intrinsic Flat Distance between
oriented Riemannian manifolds may be bounded as follows:
\be \label{bounding}
d_{\mathcal{F}}(M_1, M_2)
\le \left(\frac{n+1}{2}\right)\,\lambda^{n-1}(\lambda-1) \max\{D_1, D_2\}(V_1+A_1)\\
\ee
where  $\lambda=e^{d_{\mathcal{L}}(M_1,M_2)}$, 
$D_i=\diam(M_i)$, $V_i=\vol_n(M_i)$ and $A_i=\vol_{n-1}(\partial M_i)$.
This is proven using the theory of currents.

In work of the second author
and S.\ Lakzian, the following more constructive statement is proven
using Riemannian geometry \cite{Lakzian-Sormani}.  
Here the statement provides the
explicit Riemannian manifold $Z$ which in this case is the filling
manifold, $B$, and 
\be
A=\partial B \setminus (\psi_1(M_1)\cup \psi_2(M_2)).
\ee
{\em Suppose $M_1$ and $M_2$ are oriented precompact
Riemannian manifolds with a bi-Lipschitz map $\varphi: M_1\to M_2$
and suppose there exists $\epsilon>0$ such that
\be
(1+\varepsilon)^{-2} g_1(V,V) < g_2(\varphi_*V,\varphi_*V)
<(1+\varepsilon)^2 g_1(V,V) \qquad \forall \, V \in TM.
\ee
Then for any
\be
t_2-t_1> \frac{\arccos(1+\varepsilon)^{-1} }{\pi}\max\{\diam(M_1),\diam(M_2)\}
\ee
there is a pair of metric isometric embeddings
$\psi_i:M_i \to Z=\bar{M} \times [t_1, t_2]$ with a metric
\be \label{cor-hem-ineq}
g' \ge dt^2 +\max_{i=1,2} \cos^2((t-t_i)\pi/D_i) g_i \; \textrm { on } \; Z
\ee
and
\be\label{cor-hem-eq}
g' =  dt^2 + g_i \; \textrm  { on } \; \psi_i(M)=M \times \left\{t_i\right\}.
\ee
Thus the Intrinsic Flat distance 
between the manifolds
\be\label{77}
d_{\mathcal{F}}(M_1, M_2) \le 2|t_2-t_1|
\left(V_1+ V_2 + A_1+A_2\right),
\ee
where $D_i=\diam(M_i)$, $V_i=\vol_n(M_i)$ 
and $A_i=\vol_{n-1}(\partial M_i)$.}
In the proof of this theorem, the cosine term in (\ref{cor-hem-ineq}) 
arises from a comparison to the equator isometrically embedded
in a sphere of diameter $D_i$.

Constructive estimates on the Intrinsic Flat Distance are
also provided in previous work of the authors studying the
stability of the Positive Mass Theorem \cite{Lee-Sormani1}.
In that paper, no bi-Lipschitz maps are constructed.  Instead
metric isometric embeddings are created from
Riemannian isometric embeddings using strips whose widths
are bounded by controlling the embedding constant.    

\subsection{Intrinsic Flat Estimates}

Theorem~\ref{thm-penrose-IF} now follows from 
Proposition~\ref{Lip-Tube} combined with 
the bound on Intrinsic Flat Distance in terms of Lipschitz Distance.

\begin{proof}[Proof of Theorem~\ref{thm-penrose-IF}]
Let $r_i:=\left(\frac{A_i}{\omega_{n-1}}\right)^{\frac{1}{n-1}}$.
By (\ref{bounding}) combined with Proposition~\ref{Lip-Tube}, we know
that for any $\epsilon>0$, we can find $\delta>0$
as in Proposition~\ref{Lip-Tube} so that
\be 
d_{\mathcal{F}}\left(T_D(\Sigma_1), T_D(\bar{\Sigma}_1)\right)
\le \left(\frac{n+1}{2}\right)\,e^{(n-1)\epsilon}(e^\epsilon-1) 
\max\left\{D_1, D_2\right\}(V_1+A_1)
\ee
where by the monotonicity of the areas of the level sets, we have
\begin{eqnarray}
D_1&=&\diam\left(T_D(\Sigma_1)\right) \le 2D + \pi r_1\\
D_2&=&\diam\left(T_D(\bar{\Sigma}_1)\right) \le 2D + \pi r_1\\
V_1&=&\vol_n\left(T_D(\bar{\Sigma}_1)\right) \\
&\le & D A_0 +\vol\left(\rsch^{-1}[r_0, r_1+D])\subset \Msch(m_0)\right)\\
A_1 &=& \vol_{n-1}\left(\partial T_D(\bar{\Sigma}_1) \right) \\
& \le & 2 \vol_{n-1}\left(\rsch^{-1}(r_1+D))\subset \Msch(m_0) \right).
\end{eqnarray}
Since all these terms are bounded uniformly depending 
only on $r_0$, $r_1$ and $m_0$, which in turn depend
only on $A_0$ and $A_1$ as in the proof of 
Theorem~\ref{thm-penrose}, the result follows. 
\end{proof}

\section{Further Remarks}

In the first subsection we explain why it is necessary to introduce a
notion of {\em depth} and use \emph{appended} Schwarzschild manifolds in the statements of our main theorems.  We also explain why Theorem \ref{thm-penrose} cannot be strengthened to $C^2$ closeness of the metrics. In the second subsection, we consider the stability of the 
Positive Mass Theorem in \cite{Lee-Sormani1} and contrast 
it with the situation for the Penrose Inequality.  In the third subsection, we present 
Example~\ref{ex-many-wells} demonstrating the necessity of rotational
symmetry to obtain the Lipschitz bounds in 
Theorem~\ref{thm-penrose}, followed by
Remark~\ref{rmrk-many-wells} explaining heuristically why
these examples may be controlled using the Intrinsic Flat distance. 
We close with a subsection describing open problems and conjectures.

\subsection{Rotationally Symmetric Examples}\label{examples}

First we describe the examples depicted in Figure~\ref{fig-thm-penrose}.  Recall from Section \ref{background}
that for rotationally symmetric manifolds with outermost minimal boundaries, monotonicity of the Hawking mass 
of the symmetric spheres is \emph{equivalent} to nonnegativity of scalar curvature.  Therefore, we have the following lemma which will be useful
for constructing examples:

\begin{lem}\label{lem-ex}
{\em There is a bijection between elements of
$\RSB_n$ and increasing functions
$m_H:[r_0,\infty)\to\R$ such that 
\be \label{lem-ex-1}
m_H(r_0)=\frac{1}{2}r_0^{n-2}
\ee 
and
\be
m_H(r)<\frac{1}{2}r^{m-2}
\ee 
for $r>r_0$.  In this section we will call these functions
\emph{admissible Hawking mass functions}.  }
\end{lem}
This lemma was proven in \cite{Lee-Sormani1}.

\begin{example}\label{ex-penrose-well}
Let $n\ge3$.  Given any $0<A_0<A_1$, and any large $L>0$ and small $\delta>0$, 
there exists $M\in \RSB_n$ with $|\partial M|=A_0$
such that 
   \be    \madm \le \frac{1+\delta}{2} \left(\frac{A_0}{\omega_{n-1}}\right)^{\frac{n-2}{n-1}},   \ee 
and also
   \be    d(\partial M, \Sigma_1) > L,    \ee 
where $\Sigma_1$ is the symmetric sphere of area $A_1$.  
\end{example}

\begin{rmrk}This example shows that the appended Schwarzschild manifolds must be included in Theorem \ref{thm-penrose}.
\end{rmrk} 

\begin{proof}
Define
\begin{align}
 r_0&:=\left(\frac{A_0}{\omega_{n-1}}\right)^{\frac{1}{n-1}}\\ 
 r_1&:=\left(\frac{A_1}{\omega_{n-1}}\right)^{\frac{1}{n-1}}.
 \end{align} 
 (This is not the same $r_1$ from Section \ref{lipschitz}.) Then 
   \be    d(\partial M, \Sigma_1)  = \int_{r_0}^{r_1} \frac{ds}{dr}\,dr.   \ee 
Let $\epsilon>0$ and define $\underline{r}_\epsilon$ and $\overline{r}_\epsilon$ so that
 \begin{align}
 \underline{r}_\epsilon^{n-2}(1-\epsilon)&=r_0^{n-2}\\
 \overline{r}_\epsilon^{n-2}(1-\epsilon)&=(1+\delta)r_0^{n-2}.
 \end{align}
   We define a continuous admissible Hawking function $m_H:[r_0\,\infty)\to\infty$ as follows:
\begin{equation}
m_H(r) = \left\{ 
\begin{aligned}
\tfrac{1}{2}r_0^{n-2} &\text{for }r_0\le r \le \underline{r}_\epsilon \\
\tfrac{1}{2}r^{n-2}(1-\epsilon) &\text{for }\underline{r}_\epsilon \le r\le \overline{r}_\epsilon\\
\tfrac{1}{2}r_0^{n-2}(1+\delta) &\text{for }\overline{r}_\epsilon \le r 
\end{aligned}\right.
\end{equation}
Be Lemma \ref{lem-ex}, this admissible Hawking function defines an element of $\RSB_n$, and one can easily see that it has all of the desired properties except for $d(\partial M, \Sigma_1) > L$.  To see that we can choose $\epsilon$ small enough to make this true, observe that by the definition of $m_H$ \eqref{define-hawking}, we have $\frac{dr}{ds}=\sqrt{\epsilon}$ over the region $\underline{r}_\epsilon \le r\le \overline{r}_\epsilon$.  Assuming that $\delta$ is small enough so that $\overline{r}_\epsilon\le r$, we have
\begin{align}
 d(\partial M, \Sigma_1)&  \ge \int_{\underline{r}_\epsilon}^{\overline{r}_\epsilon} \frac{1}{\sqrt{\epsilon}}\,dr\\
  &=  \frac{r_0}{\sqrt{\epsilon}} \frac{ (1+\delta)^{\frac{1}{n-2}}-1}{(1-\epsilon)^{\frac{1}{n-2}}},
 \end{align}
which tends to infinity as $\epsilon\to0$.  Note that this example can be modified in a straightforward manner to obtain a smooth example.  
\end{proof}

The next example proves that
the Lipschitz convergence cannot be improved to $C^2$
convergence in Theorem~\ref{thm-penrose}.
\begin{ex} \label{ex-not-C2}
There exists a sequence $M_j\in \RSB_n$ such that each $|\partial M_j| = A_0$, 
   \be    \lim_{j\to\infty} \madm(M_j) =  \left(\frac{A_0}{\omega_{n-1}}\right)^{\frac{n-2}{n-1}},   \ee 
and $M_j$ converges to $\Msch(m_0)$ in the Lipschitz sense, but not in the $C^2$ sense, where 
$m_0=\left(\frac{A_0}{\omega_{n-1}}\right)^{\frac{n-2}{n-1}}$.
\end{ex}
\begin{proof}
Choose a sequence $m_j\searrow m_0$.  For each $j$ we construct an admissible Hawking function $m_H^{(j)}$ such that $m_H^{(j)}(r)=m_0$ for $r\in[r_0, 2r_0]$, and then $m_H^{(j)}(r)$ takes a a very brief sharp turn upward so that $\frac{m_H^{(j)}}{dr}  > \frac{1}{j}$ for some $r\in(2r_0,3r_0)$, and then $m_H^{(j)}=m_j$ for $r\in[3r_0,\infty)$.  By our proof Theorem \ref{thm-penrose}, the corresponding manifolds must converge to $\Msch(m_0)$ in the Lipschitz sense.  (We know that the depths are zero, because each $M_j$ is exactly Schwarzschild in a fixed neighborhood of $\partial M_j$.)  But \eqref{define-hawking} tells us that
   \be   \frac{m_H^{(j)}}{dr} = \frac{r^{n-2}}{2}\mathrm{R},   \ee 
and therefore $\mathrm{R}$ can be arbitrarily large as $j\to\infty$, so $C^2$ convergence is impossible.
\end{proof}

\subsection{Contrasting with Positive Mass Stability}

Theorem~\ref{thm-penrose-IF} is not a stability statement,
because the appended Schwarzschild space depends on the depth of $M$.  In particular, this means that if we have a sequence of manifolds $M_j$ with fixed $|\partial M_j|$ whose $\delta_j$'s approach zero, the sequence need not converge to anything, since their depths need not converge (which is intuitively clear from \ref{ex-penrose-well}).  In this sense Theorem~\ref{thm-penrose-IF} is weaker than the Positive Mass Stability Theorem proven by the authors in \cite{Lee-Sormani1}.   There we
proved that if $M \in \RSB_n$ has ADM mass close to zero, then tubular
neighborhoods in $M$ are close in the Intrinsic Flat sense to 
tubular neighborhoods in Euclidean space.   In particular, a sequence $M_j$ whose masses approach zero must converge to Euclidean space in the pointed Intrinsic Flat, with appropriately chosen basepoints.

On the other hand Theorem~\ref{thm-penrose} 
is far stronger than the Positive
Mass Stability Theorem in \cite{Lee-Sormani1} in the sense that we are able
to obtain Lipschitz estimates rather than just Intrinsic Flat estimates.
The reason one cannot obtain Lipschitz estimates for the near-equality case of 
the Positive Mass Theorem is that there $\partial M$ can be both arbitrarily small and arbitrarily deep.
If we look at a limit of such examples with masses approaching zero, we see that the deep regions are not becoming close to a cylinder but rather a line segment.  There cannot
be a bi-Lipschitz map from a region of dimension $\ge 3$, to
a region of dimension $1$.   However, this interior region can be controlled
by constructing a filling manifold and estimating its volume, leading to a stability result
formulated in terms of Intrinsic Flat Distance \cite{Lee-Sormani1}.

Note that in \cite{Lee-Sormani1}, it was possible to obtain a Lipschitz estimate \emph{away} from the thin, deep well, which could have been proved in a manner similar to Section \ref{lipschitz}.  Such a proof could have been used as a step in proving the main results of \cite{Lee-Sormani1}, but instead we chose a proof that emphasized features of Intrinsic Flat Distance which are better suited to more general problems.   
Even in the Penrose setting, the Lipschitz estimate in Theorem \ref{thm-penrose} is special to the rotationally
symmetric case, as seen below.

\subsection{Many Wells}

In this subsection we present 
Example~\ref{ex-many-wells} demonstrating the necessity of rotational
symmetry in Theorem~\ref{thm-penrose}, followed by
Remark~\ref{rmrk-many-wells} explaining heuristically why
these example with many wells can be controlled
with the Intrinsic Flat distance.   

\begin{example}\label{ex-many-wells}
There exist sequences of
Riemannian manifolds satisfying all the hypotheses of
Theorem~\ref{thm-penrose} except the rotational symmetry,
with $\delta_j \to 0$, such that the sequence does not become Lipschitz close (or even Gromov-Hausdorff close) to 
$\Msch(m_0, L)$ for any $L$.  
\end{example}

\begin{proof}
One may construct a Riemannian manifold with many wells similar
to the example constructed in \cite{Lee-Sormani1}.  As in that example,
one begins with a rotationally symmetric manifold with thin regions
of constant sectional curvature.  The only difference is that here
we start with manifolds $M'_j\in \RSB_n$ that have fixed 
boundary of (hyper-)area $A_0$
and edit in thin strips of constant curvature while ensuring that
\be
\madm(M'_j) \to m_0= \frac{1}{2} \left(\frac{A_0}{\omega_{n-1}}\right)^{\frac{n-2}{n-1}}.
\ee  
We can do this in a way such that $M'_j$ converge in the Lipschitz
sense to some $\Msch(m_0)$.   

Next we edit in rotationally symmetric wells of some fixed depth $h>0$ replacing small balls in the strips of constant sectional curvature.  We can
edit in arbitrarily many such wells just as in the example constructed in
\cite{Lee-Sormani1}.  This creates new manifolds $M_j$ which
satisfy all the conditions of Theorem~\ref{thm-penrose} except the
rotational symmetry.  As they have increasingly many wells of
depth $h$, they have no Gromov-Hausdorff limit and do not 
converge in the Lipschitz sense.
\end{proof}

\begin{rmrk}\label{sum-Hawking} 
Suppose that $M^n$ is asymptotically flat with
positive scalar curvature and that $\partial M$ is outermost
minimizing.  Suppose further
that $M$ contains many rotationally symmetric wells, $W_i$, as in Example~\ref{ex-many-wells}.  If $n<8$,  the Penrose inequality guarantees that 
\be
\sum_i m_H(\partial W_i)^{\frac{n-1}{n-2}} \le \madm(M)^{\frac{n-1}{n-2}} - m_H(\partial M)^{\frac{n-1}{n-2}}
\ee
To see this, replace each well $W_i$ by a region of  $\msch(m_H(W_i))$ to create a new manifold $N$, and then apply the Penrose inequality to $N$ \cite{Bray-Lee}.
\end{rmrk}

\begin{rmrk}\label{rmrk-many-wells}
Example~\ref{ex-many-wells} does not appear to contradict an
extension of Theorem~\ref{thm-penrose-IF} to the setting
where there is no rotational symmetry.  
One can imagine constructing fillings for each well and  
controlling the volume of each filling in terms of the
Hawking mass at the top of each well
using techniques similar to those developed in \cite{Lee-Sormani1}.  
\end{rmrk}

\subsection{Conjectures}

While Theorem \ref{thm-penrose} does not generalize, we hope that Theorem \ref{thm-penrose-IF} does.  We propose a conjecture along the lines of what we conjectured for stability of the Positive Mass Theorem in \cite{Lee-Sormani1}. 

\begin{defn}
Let $\mathcal{M}$ be a subclass of asymptotically flat
three dimensional Riemannian manifolds of nonnegative
scalar curvature, whose boundaries are outermost minimal surfaces.
\end{defn}

\begin{conj}\label{conj-penrose-IF}
Let $0<A_0<A_1$.   For any $D>0$ and $\epsilon>0$, there exists $\delta>0$ depending on these constants, such that for any $M^3\in\mathcal{M}$ satisfying
   \be    
   \madm \le (1+\delta)\sqrt{\frac{A_0}{16\pi}},   
   \ee 
   where
$A_0$ is the (hyper-)area of one connected component of $\partial M$, 
then
   \be    
   d_{\mathcal{F}}\left(T_D(\Sigma_1)\subset M, T_D(\bar{\Sigma}_1)\subset \Msch^n (m_0, L)\right) < \epsilon,  
    \ee 
where  $\Sigma_1$ is a \emph{special surface} of area $A_1$, 
$L$ is the \emph{depth} of $M$, 
$\bar{\Sigma}_1$ is the symmetric sphere of area $A_1$ in $\Msch(m_0,L)$
and $m_0= \frac{1}{2} \left(\frac{A_0}{\omega_{n-1}}\right)^{\frac{n-2}{n-1}}$. 
\end{conj}

We are deliberately vague as to how restrictive the class $\mathcal{M}$ needs to be.  The conjecture may
require uniform conditions at infinity.
We have also been vague as to what the {\em special surface}, $\Sigma_1$,
should be.   The main point about $\Sigma_1$ is that it should avoid the wells and also not escape to infinity.  
Please see \cite{Lee-Sormani1} for more discussion about how $\Sigma_1$ might be defined.   We have also been vague as to what
the appropriate definition of depth should be in the general setting.
See Remark~\ref{rmrk-depth} for further discussion regarding this notion of
depth.

\bibliographystyle{plain}
\bibliography{2011}

\end{document}